\documentclass[12pt]{article}
\usepackage{amsmath}
\usepackage{amssymb,color}

\def\0{{\bf 0}}
\def\1{{\bf 1}}
\def\l{{\ell}}

\def\proof{\noindent{\bf Proof: }}
\def\qed{ \hskip 20pt{\vrule height7pt width6pt depth0pt}\hfil}
\def\forb{{\mathrm{forb}}}
\def\ext{{\mathrm{Ext}}}
\def\Av{{\mathrm{Avoid}}}

\newcommand{\linelessfrac}[2]{\genfrac{}{}{0pt}{}{#1}{#2}}

\newcommand{\ncols}[1]{\| #1 \|}

\newcommand{\rf}[1]{(\ref{#1})}
\newcommand{\trf}[1]{Theorem~\ref{#1}}

\newcommand{\lrf}[1]{Lemma~\ref{#1}}

\newcommand{\srf}[1]{Section~\ref{#1}}

\newtheorem{thm}{Theorem}[section]
\newtheorem{lemma}[thm]{Lemma}

\newtheorem{conj}[thm]{Conjecture}
\newtheorem{defin}[thm]{Definition}
\newenvironment{E}{\begin{equation}}{\end{equation}}
\usepackage[english]{babel}
\makeatother

\usepackage{graphicx} 
\title{
   Exact Bounds for Forbidden Configurations and the Extremal Matrices}
\author{R.P. Anstee\thanks{Research supported in part by
NSERC}, Oakley Edens\thanks{Research supported in part by NSERC USRA}, Arvin Sahami\thanks{Research supported in part by NSERC USRA}, Jaehwan Seok\thanks{Research supported in part by NSERC USRA} \\ Mathematics Department\\The University of British Columbia\\Vancouver,
B.C. Canada V6T 1Z2\\ {\small{\texttt{anstee@math.ubc.ca}}}, 
{\small{\texttt{oedens@math.harvard.edu}}},{\small{\texttt{sahamiarvin52@gmail.com}}},\\
{\small{\texttt{jseok627@student.ubc.ca}}} 
\and Attila Sali\thanks{Research  partially supported by the
    National Research, Development and Innovation Office (NKFIH)
    grants K--132696 and SNN-135643. 
}\\ HUN-REN Alfr\'ed R\' enyi Institute of Mathematics\\Budapest, Hungary and \\ Department of Computer Science\\ Budapest University of Technology and Economics\\ {\small{\texttt{sali.attila@renyi.hu}}} 
} 

\begin{document}

\maketitle

\begin{abstract}
    Let $F$ be a $k\times \l$ (0,1)-matrix. A matrix is {\it simple} if it is a (0,1)-matrix with no repeated columns. A (0,1)-matrix $A$ to said to have a $F$ as a \emph{configuration} if there is a submatrix of $A$ which is a row and column permutation of $F$. In the language of sets, a configuration is a \emph{trace}.
 Let $\Av(m,F)$ be all simple $m$-rowed matrices $A$ with no configuration $F$.  Define $\forb(m,F)$ as the maximum number of columns of any matrix in $\Av(m,F)$.   The $2\times (p+1)$ (0,1)-matrix $F(0,p,1,0)$ consists of a row of $p$ 1's and a row of one 1 in the remaining column. The paper determines  $\forb(m,F(0,p,1,0))$ for $1\le p\le 9$ and the extremal matrices are characterized. A construction may be extremal for all $p$.

\end{abstract}

\section{Introduction}

The paper considers bounds for $F(0,p,1,0)$ where
\begin{E}F(0,p,1,0)=\left[\begin{array}{c@{}}\\ \\ \end{array}\right. 
\overbrace{
\begin{array}{@{}c@{}}1\,\,1\cdots 1\,\\ 0\,\,0\cdots 0\,\\ \end{array}}^p 
\begin{array}{@{}c@{}}0\\ 1\\\end{array} 
\left.\begin{array}{@{}c}\\ \\ \end{array}\right].
\label{F(r,p,q,s)}\end{E}

An $m\times n$ matrix $A$ is said to be \emph{simple} if it is a (0,1)-matrix with no repeated columns. There is a natural correspondence between columns of $A$ and subsets of $[m]$. We use matrix terminology as follows. Let $\ncols{A}$ be the number of columns of $A$. For a given matrix $F$, we say $F$ is a \emph{configuration} in $A$, denoted $F\prec A$, if there is a submatrix of $A$ which is a row and column permutation of $F$. Define 
$$\Av(m,F)=\left\{A\,|\,A \hbox{ is }m\hbox{-rowed and simple}, F\not\prec A\right\},$$ 
$$\forb(m,F)=\max_{A\in\Av(m,F)}\ncols{A}.$$
A matrix $A\in\Av(m,F)$ is called {\emph{extremal}} if $\ncols{A}=\forb(m,F)$ and let 
$$\ext(m,F)=\{A\in\Av(m,F)\,|\,\ncols{A}=\forb(m,F)\}.$$ 
Many results are in a survey \cite{survey} including for 2-rowed $F$.

For a subset $S$ of rows, define $A_S$ as the submatrix of $A$ formed by those rows.   Some important matrices include  $K_k$,  the $k\times 2^k$ matrix of all possible (0,1)-columns on $k$ rows. We denote by $F^c$ the (0,1)-complement of $F$ so that $I_k^c$ is the complement of the identity. Define $K_k^s$ to be the $k\times \binom{k}{s}$ matrix of columns of sum  $s$. We define $\1_k$ as $k\times 1$ column of 1's, $\0_k$ as $k\times 1$ column of 0's. 

Problems in extremal combinatorics  are first concerned with bounds but considerations of what happens at the bound are also explored.  
The paper determines bounds $\forb(m,F(0,p,1,0))$ for $3\le p\le 9$ and provides a characterization of  matrices in  $\ext(m,F)$  (\srf{F2(0,3,1,0)} and \srf{F2(0,6,1,0)}). General lemmas about $\Av(m,F(0,p,1,0))$ are in \srf{F(0,p,1,0)}.  Construction \rf{construction} seems crucial for determining $\ext(m,F)$.  
  Graph theory with directed and undirected edges is used to describe  a matrix  in $\Av(m,F)$ \rf{wim}. \srf{F2(0,p,1,0)} and \srf{cliquecomponents} have proofs of the important Lemmas. The most important result is \lrf{lemma:cliques} which asserts that the components induced by undirected edges are cliques of undirected edges. This yields bounds and characterizations of $\ext(m,F)$.   The arguments for  \lrf{lemma:cliques} get more complicated as $p$ grows.  We ignore row and column permutations of our matrices unless explicitly stated.

A detailed characterization of $\ext(m,F(0,3,1,0))$ is in \srf{F2(0,3,1,0)}. Bounds for $F(0,4,1,0)$, $F(0,5,1,0)$ are in \cite{survey}.
 We obtain  new exact bounds   in \srf{F2(0,6,1,0)} for $p\in\{6,7,8,9\}$ which also yield characterization of \hfil\break$\ext(m,F(0,p,1,0))$ most notably in \lrf{extremalclique}. 
\begin{thm} $\forb(m,F(0,6,1,0))\le \lfloor\frac{21}{5}m\rfloor+1$ with equality only for $m\equiv 0(\hbox{mod }5)$. Also $\forb(m,F(0,6,1,0))= \lfloor\frac{21}{5}m\rfloor$ for $m\equiv 1(\hbox{mod }5)$ and $m\ge 6$. \label{F61}\end{thm}
We characterize matrices in \hfil\break $\ext(m,F(0,6,1,0))$ for $m\equiv 0,1 (\hbox{mod }5)$. 

\begin{thm} $\forb(m,F(0,7,1,0))\le \lfloor\frac{24}{5}m\rfloor+1$ with equality for $m\equiv 0(\hbox{mod }5)$. \label{F71}\end{thm}

\begin{thm} $\forb(m,F(0,8,1,0))\le \lfloor\frac{27}{5}m\rfloor+1$ with equality for $m\equiv 0(\hbox{mod }5)$. \label{F81}\end{thm}

\begin{thm} $\forb(m,F(0,9,1,0))\le \lfloor\frac{31}{5}m\rfloor+1$ with equality for $m\equiv 0(\hbox{mod }5)$. \label{F91}\end{thm}

\section{Ideas for $\Av(m,F(0,p,1,0))$}\label{F(0,p,1,0)}

We follow the proof ideas of \cite{AKa}. For $A\in\Av(m,F(0,p,1,0))$, then on any pair of rows $i,j$ we  have one of three cases:
\begin{E} \left(\begin{array}{@{}c@{}}\\ \left.\begin{array}{c@{}}i\\ j\\ \end{array}\right.\\ \end{array}
\begin{array}{@{}c}\le p-1\\\left[\begin{array}{c}0\\1\\ \end{array}\right]\\ \end{array}
\hbox{and}\begin{array}{c@{}}\\ \left.\begin{array}{c@{}}i\\ j\\ \end{array}\right.\\ \end{array}
\begin{array}{@{}c}\le p-1\\\left[\begin{array}{c}1\\0\\ \end{array}\right]\\ \end{array}\right)\,\,
\hbox{ or }\begin{array}{c@{}}\\ \left.\begin{array}{c@{}}i\\ j\\ \end{array}\right.\\ \end{array}
\begin{array}{@{}c}\hbox{no}\\\left[\begin{array}{c}0\\1\\ \end{array}\right]\\ \end{array}
\hbox{ or }\begin{array}{c@{}}\\ \left.\begin{array}{c@{}}i\\ j\\ \end{array}\right.\\ \end{array}
\begin{array}{@{}c}\hbox{no}\\\left[\begin{array}{c}1\\0\\ \end{array}\right]\\ \end{array}.
\label{wim}\end{E}
This is the `What is Missing' idea in \cite{survey}. Note the row ordering.  We form a graph $G(A)$ on rows of $A$ of edges and directed edges as follows: 
\begin{itemize}
    \item undirected edge $i-j$ if there are at most $p-1$ 
    $\linelessfrac{i}{j}\left[\linelessfrac{0}{1}\right]$ and at most $p-1$ $\linelessfrac{i}{j}\left[\linelessfrac{1}{0}\right]$ submatrices on  rows $\{i,j\}$
   \item directed edge $i\rightarrow j$ if there is no submatrix 
   $\linelessfrac{i}{j}\left[\linelessfrac{0}{1}\right]$ on  rows $\{i,j\}$ \end{itemize} 
   Every pair of rows of $G(A)$ must either be joined by a directed or undirected edge   to avoid $F(0,p,1,0)$.  Assume that if we have both $i-j$ and $i\rightarrow j$, we will ignore $i-j$.

   Let the components induced by  the undirected edges be $C_1,C_2,\ldots, C_t$. 
 This yields  a promising construction for
   an $A\in\Av(m,F(0,p,1,0))$. The entry $\0$  denotes a block or column of 0's and $\1$ denotes a block or column of 1's.  Let $B_i$ be the submatrix,  on  the  rows of the component $C_i$, that has the non-constant columns on rows $C_i$ plus perhaps the column of 0's.

   \begin{lemma}
   Let $B_i$ be a $k_i$-rowed simple matrix with $B_i \in \Av(k_i,F(0,p,1,0))$ and assume $B_i$ has  no column of 1's but may have a column of 0's. Let
  \begin{E}A=
\begin{array}{c@{}}
C_1\\ C_2\\ C_3\\ \vdots\\ \vdots\\ C_t\\ \end{array}\left[\begin{array}{cccccccc}
B_1&{\1}&{\1}&&\cdots&&{\1}&{\1}\\
{\0}&B_2&{\1}&&\cdots&&{\1}&\1\\
{\0}&{\0}&B_3&&\cdots&&{\1}&\1\\
\vdots&\vdots&\vdots&&&\ddots&&\vdots\\
\vdots&\vdots&\vdots&&&&\ddots&{\1}\\
{\0}&{\0}&&\cdots&&{\1}&B_t&\1\\
\end{array}\right]\label{construction}\end{E} 
    Then  $A\in\Av(m,F(0,p,1,0))$ with $m=\sum_{i=1}^tk_i$ and $\ncols{A}=1+\sum_{i=1}^t\ncols{B_i}$.    
   \label{mainconstruction}\end{lemma}
\proof  Taking rows $r\in C_i$, $s\in C_j$ for $i<j$ we have no  submatrix $\linelessfrac{r}{s}\left[\linelessfrac{0}{1}\right]$. Since $B_i \in \Av(k_i,F(0,p,1,0))$, we deduce that $A\in\Av(m,F(0,p,1,0))$. \qed
   \vskip 10pt
If the directed edges of $G(A)$ have transitivity  for $A\in\ext(m,F(0,p,1,0))$, then there is an ordering of the components (induced by the undirected edges) so that the directed edges go from $B_i$ to $B_j$ for $i<j$.  This yields the construction \rf{construction}.  We would like to show that in $A$  the components formed by the undirected edges are cliques.  \lrf{lemma:cliques} handles this.
Determining the best choices for $B_i$ in \rf{construction} is somewhat like a design problem as seen in \lrf{risefall} Upper Bound Lemma.  The best choices seem to have the rows of $B_i$ being a clique in $G(A)$ as proven in \lrf{lemma:cliques} yielding a bound $\forb(m,F(0,p,1,0))\le \lfloor c_pm\rfloor +1$ with $c_p$ given in \rf{cp}.

For $p=2$, we can take $B_i=[K_2\backslash \1_2]$ yielding  $\frac{3}{2}m+1$ columns for $m\equiv 0({\hbox{mod }}2)$ which  easily yields $\lfloor\frac{3}{2}m\rfloor+1$. For $p=3$, take $B_i=[K_3^2K_3^1K_3^0]$ yielding $\frac{7}{3}m+1$ columns for $m\equiv 0({\hbox{mod }}3)$. For $p=4$, take $B_i=[K_4^3DK_4^1K_4^0]$ where $D$ is the $4\times 2$ matrix of 2 complementary columns of sum 2 yielding $\frac{11}{4}m+1$ columns for $m\equiv 0({\hbox{mod }}4)$. For $p=5$, take $B_i=[K_4^3K_4^2K_4^1K_4^0]$ yielding $\frac{15}{4}m+1$ columns for $m\equiv 0({\hbox{mod }}4)$.   For $p=6$, take $B_i=[K_5^4K_5^2K_5^1K_5^0]$ in \rf{construction} yielding $\frac{21}{5}m+1$ columns for $m\equiv 0({\hbox{mod }}5)$.  We could also take  $[K_5^4K_5^3K_5^1K_5^0]$ (essentially the complement) or indeed $[K_5^4G_5(G_5')^cK_5^0]$ where $G_5$ is the incidence matrix of a cycle on 5 vertices and $G_5'$ is the incidence of the cycle that is the graph complement of $G_5$.   These three possibilities $K_5^2$, $K_5^3$ and $[G_5(G_5')^c]$ have exactly 3 configurations $\linelessfrac{i}{j}\left[\linelessfrac{0}{1}\right]$ for each possible ordered pair of rows $i,j$. \trf{F61} shows 
$\forb(m,F(0,6,1,0))=  \frac{21}{5}m+1$ for  $m\equiv 0(\hbox{mod }5)$ and hence for $m\equiv 0(\hbox{mod }5)$ we have constructions for $\ext(m,F(0,6,1,0))$ with three choices for each $5\times 21$ $B_i$ and moreover it is the only way to construct an $A\in\ext(m,F(0,6,1,0))$.

For $p=7$, a computer search yields $\forb(5,F(0,7,1,0))=25$ with a construction:
\begin{E}\left[K_5^0\,K_5^1\,\,\begin{array}{ccccccccccccc}
1&1&1&1&0&0&0&0&0&1&1&0&0\\
1&0&0&0&1&1&1&0&0&1&0&1&0\\
0&1&0&0&1&0&0&1&1&1&0&0&1\\
0&0&1&0&0&1&0&1&0&0&1&1&1\\
0&0&0&1&0&0&1&0&1&0&1&1&1\\
\end{array}\,\,K_5^4\,K_5^5\right].\label{F7construction}\end{E}
This meets the bound of the Upper Bound Lemma for a clique.  Also it yields 
$\forb(m,F(0,7,1,0))\ge \frac{24}{5}m+1$ for $m\equiv 0(\hbox{mod }5)$ and hence our guess   $c_7=\frac{24}{5}$.  

For  $p=8$,
a computer search yields $\forb(5,F(0,8,1,0))=28$ with a construction:
\begin{E}\left[K_5^0\,K_5^1\,\,\begin{array}{ccccccccccccccccc}
1&1&1&1&0&0&0&0&  1&1&1&1&0&0&0&0\\
1&0&0&0&1&1&0&0&  1&1&0&0&1&1&1&0\\
0&1&0&0&1&0&1&0&  1&0&1&0&1&1&0&1\\
0&0&1&0&0&1&0&1&  0&0&1&1&1&0&1&1\\
0&0&0&1&0&0&1&1&  0&1&0&1&0&1&1&1\\
\end{array}\,\,K_5^4\,K_5^5\right].\label{F8construction}\end{E}
This meets the bound of the Upper Bound Lemma for a clique.  Also it yields 
$\forb(m,F(0,8,1,0))\ge \frac{27}{5}m+1$ for $m\equiv 0(\hbox{mod }5)$ and hence our guess  $c_8=\frac{27}{5}$.

For $p=9$,  
the construction $K_5$ yields $\forb(5,F(0,9,1,0))=32$. This meets the bound of the Upper Bound Lemma for a clique. Also it yields 
$\forb(m,F(0,9,1,0))\ge \frac{31}{5}m+1$ for $m\equiv 0(\hbox{mod }5)$. 
and hence our guess  $c_9=\frac{31}{5}$.

From our constructions we define
  
 \begin{E}
\hbox{\hfil
\vbox{\offinterlineskip
\hrule
\halign{&\vrule#&
\strut\,\,\hfil#\hfil\,\,\cr
height2pt&\omit&&\omit&&\omit&&\omit&&\omit&&\omit&&\omit&\cr
&$p=3$&&$p=4$&&$p=5$&&$p=6$&&$p=7$&&$p=8$&&$p=9$&\cr
height2pt&\omit&&\omit&&\omit&&\omit&&\omit&&\omit&&\omit&\cr
\noalign{\hrule}
height2pt&\omit&&\omit&&\omit&&\omit&&\omit&&\omit&&\omit&\cr
&$c_3=\frac{7}{3}$&&$c_4=\frac{11}{4}$&&$c_5=\frac{15}{4}$&&$c_6=\frac{21}{5}$&&$c_7=\frac{24}{5}$&&$c_8=\frac{27}{5}$&&$c_9=\frac{31}{5}$&\cr
height2pt&\omit&&\omit&&\omit&&\omit&&\omit&&\omit&&\omit&\cr
\noalign{\hrule}
}
}
}  
 \label{cp}\end{E}

The matrix $K_t$ meets the bound of the Upper Bound Lemma \lrf{risefall} for a clique on $2^{t-2}+1$ vertices  and  $B_i=K_t\backslash \1_t$  can be used in \rf{construction}.  We make the following conjecture:

\begin{conj}Let $t$ be given. $\forb(m,F(0,2^{t-2}+1,1,0))\le \frac{2^{t}-1}{t}m+1$ with equality for $m\equiv 0(\hbox{mod }t)$.\end{conj}

\section{Lemmas for $F(0,p,1,0)$}\label{F2(0,p,1,0)}

We will use Lemmas here and in \srf{cliquecomponents} to show that $\forb(m,F(0,p,1,0))\le c_pm+1$
for $p=3,4,\ldots,9$.   To show $\forb(m,F(0,p,1,0))=\lfloor c_p m\rfloor +1$  then $\ncols{B_i}$ needs to be close to 
   $c_p|B_i|$.  The following is helpful in inductive proofs.

   \begin{defin} Let $A\in\Av(m,F(0,p,1,0))$. If we delete $s$ rows and $t$ columns and the resulting matrix $A_{m-s}$ is simple then we assign a cost to the deletion
$$\hbox{cost}=c_p s\,\,-t\,.$$
Given a component $C_i$ of $G(A)$ induced by undirected edges, let $B_i$ be the matrix formed from the non-constant columns of $A_{C_i}$  plus the column of 0's if present.   We can delete the $s=|C_i|$ rows $C_i$ and the $\ncols{B_i}$ columns of $B_i$ to obtain a simple matrix thus 
\begin{equation}\hbox{cost of component }C_i=c_p |C_i|-\ncols{B_i}.
\label{costcomponent}\end{equation}
 \end{defin}

The proof of $\forb(m,F(0,p,1,0))\le\lfloor c_p m\rfloor +1$ proceeds by induction to show $\forb(m,F(0,p,1,0))\le c_p m +1$. The idea is that when applying induction we will be able to prove $\forb(m,F(0,p,1,0))\le c_pm+1$ if the costs are positive or zero. The proof can ignore components whose cost is at least 1 since then we are showing (by induction) that  $\ncols{A}\le(c_p m+1)-1$, hence $\ncols{A}\le\lfloor c_p m\rfloor$.

\begin{lemma}\label{rows}{\bf Deletion Lemma}. Assume we are trying to show by induction on $m$ that    $\forb(m,F(0,p,1,0))\le c_pm+1$. Then we are done if there is a deletion of cost at least 0. 
\end{lemma}

\proof  Assume $A\in\Av(m,F)$.     Assume we have a deletion of $k$   rows of $A$ and $t$ columns  so that the resulting matrix $A_{m-k}$ is simple: $A_{m-k}\in\Av(m-k,F)$. If the cost of the deletion is at least 0 then $c_pk>t$.  Then by induction $\ncols{A_{m-k}}\le c_p(m-k)+1$ so
$\ncols{A}\le \ncols{A'}+c_pk\le c_p(m-k)+1+c_pk=c_pm+1.$ \qed
\vskip 10pt

If there is a deletion or component of cost 1 or greater then  by induction the matrix would not be extremal.   The analysis of extremal matrices can assume no  components whose deletion cost is at least 1 are present. 

The following is used extensively particularly for a clique in $G(A)$ induced by the undirected edges.

\begin{lemma} \label{risefall} {\bf Upper Bound Lemma}. Let $B$ be a $k\times n$ (0,1)-matrix with no column of all 1's and no column of all 0's.
Assume $B$ has no pair of rows which differ in more than $t$ columns i.e. $B$ has at most $t$ disjoint configurations 
$\left[\linelessfrac{0}{1}\right]$ on the same pair of rows. Then
\begin{E}\label{risefallnotsimple} n\le \frac{tk}{2}.\end{E}
If $B$ is simple and $t\ge 4$ then
\begin{E}\label{risefallsimple} n\le \Bigl\lfloor 2k+\frac{(t-4)k(k-1)}{4(k-2)}\Bigr\rfloor.\end{E}
\end{lemma}

\proof  Each column contributes at least $k-1$ 
configurations $\left[\linelessfrac{0}{1}\right]$
 in the $k$ rows.  More than $tk/2$ of
such columns would give more than $tk(k-1)/2$ of the $\left[\linelessfrac{0}{1}\right]$ configurations in $k(k-1)/2$
pairs of rows in $B$. One pair of rows would then contain more
than $t$ configurations, a contradiction yielding (\ref{risefallnotsimple}).

If $B$ is simple then we note there are at most $2k$ columns which each have only $k-1$ configurations, namely the columns with at most one 1 and the columns with at most one 0. All other columns have at least $2(k-2)$ configurations
$\left[\linelessfrac{0}{1}\right]$ with equality for column sums 2 or $k-2$. We deduce
$2k(k-1)+(n-2k)2(k-2)\le tk(k-1)/2$ and  obtain (\ref{risefallsimple}).

\qed
\vskip 10pt

This lemma applied to $B$ simple would need to be adapted for $t> 2(\binom{k}{2}+\binom{k}{1})$ 
 when the columns of sum 2 and $k-2$ have already been used up and consider the columns of sum 3, $k-3$, etc.   This would be needed for $p>9$. 

 Let $D(A)$ denote the subgraph of $G(A)$ of the directed edges.

\begin{lemma}\label{graph} {\bf Transitivity Lemma}. Let $A\in\Av(m,F(0,p,1,0))$.  Assume $2c_p\ge p$ which is true for $p=3, 4,\ldots, 9$. Each pair of rows is connected by exactly one edge of $G(A)$ else there is a deletion of cost at least 1. If $D(A)$ is not transitive then there is a deletion of cost at least 1.   Thus we may assume  that  the graph of the  directed edges $D(A)$ is transitive.
\end{lemma}

\proof It is clear in \rf{wim} that each pair $i,j$ is joined by some edge: 
$i\rightarrow j$, 
$i- j$, or
$j\rightarrow i$. Our definition of $i- j$ ensures that we do not have $i\rightarrow j$ or $j\rightarrow i$.  If  $i\rightarrow j$ and $j\rightarrow i$ 
then we can delete row $i$ and the result remains simple. The Deletion Lemma~\ref{rows} can be used with a cost of $c_p\ge 1$.

To show $D(A)$ is transitive and
contains no cycles consider the case: 
 $i\rightarrow j$ and $j\rightarrow k$. We have the three
possibilities:
$$(a)\quad \begin{matrix} i \\ \end{matrix} \begin{matrix} \nearrow \\ 
\backslash \\\end{matrix}
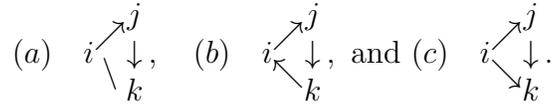
\begin{matrix} j \\ \downarrow \\ k \end{matrix} \,,\quad (b)\quad \begin{matrix} i \\
\end{matrix} \begin{matrix} \nearrow \\  \nwarrow \end{matrix} \begin{matrix} j \\
\downarrow \\ k \end{matrix} \,, \hbox{ and } (c)\quad \begin{matrix} i \\
\end{matrix} \begin{matrix} \nearrow \\  \searrow \end{matrix} \begin{matrix} j \\
\downarrow \\ k \end{matrix} \,.$$ 
For cases (a) and (b) we look
at the possible entries for these three rows. The entries above
rows $i,j,k$  indicate the number of possible columns of these types.
\begin{center}
$\begin{matrix} \\i\\j\\k \end{matrix} \;\; 
\begin{matrix} \hbox{none }\,\,\hbox{ none }\\0\cdots 0\; 0\cdots 0 \\ \;1\cdots 1\;1\cdots1
\\ 0\cdots 0\;1\cdots 1
\end{matrix}
 \;\begin{matrix} \hbox{none }\,\,\hbox{ none }\\0\cdots 0\; 1\cdots 1\\0\cdots
0\;0\cdots 0\\1\cdots 1\;1\cdots 1
\end{matrix}
 \;\begin{matrix}s\\ 1\cdots 1\\ 0\cdots 0\\ 0\cdots 0\\ \end{matrix} 
 \;\begin{matrix}t\\ 1\cdots 1\\ 1\cdots 1\\ 0\cdots 0\\ \end{matrix}
\;\begin{matrix} \\ 0\cdots 0\; 1\cdots 1\\0\cdots 0\;1\cdots 1\\0\cdots 0\;1\cdots 1
\end{matrix}$ .
\end{center}

In case (a), $s+t\le p-1$. We can delete two rows $i$ and $j$ and $s+t$ columns to obtain a simple matrix. We are done by the Deletion Lemma~\ref{rows} with a cost of $c_p\cdot 2 -(p-1)$ and hence a cost at least 1.
In case (b), rows $i,j,k$ must be identical with $s+t=0$ hence we can eliminate the two rows $i$ and $j$ and no columns to produce
a simple matrix. We are done by the Deletion Lemma~\ref{rows} with a cost of $c_p\cdot 2$ and hence a cost at least 1.
Thus we may assume $A$ can  have (c) only. Thus  $i\rightarrow j$ and $j\rightarrow k$ implies $i\rightarrow k$.    \qed
\vskip 10pt

Given that there is at most one directed edge between two rows, then Transitivity shows that there are no directed cycles.   Also we can apply  
Transitivity Lemma to $A\in\Av(m,F(0,p,1,0))$ and deduce that the components  (induced by the undirected edges) can be transitively ordered and form the structure of \rf{construction}. 

\begin{lemma}
  Let $p\in\{3,4,\ldots ,9\}$. Let $A\in\Av(m,F(0,p,1,0))$ such that each component induced  by the undirected edges is a clique. Then $\ncols{A}\le c_pm+1.$
\label{cliquesgivebound}\end{lemma}

\proof This is just the application of the Upper Bound Lemma on each component. The value of $c_p$ in \rf{cp} was determined by maximizing the function $n/k$ with $n$ and $k$ from \rf{risefallsimple}. 
\qed

\section{Clique Components}\label{cliquecomponents}

The following lemma (with a lengthy proof)  uses \lrf{cliquesgivebound} and the Upper Bound Lemma  to obtain our desired bounds. Establishing that components are cliques for $A\in\ext(m,F(0,p,1,0))$ gives a useful characterization.   In some sense the Lemma  reduces the problem to a `finite' computation for fixed $p$ since the Upper Bound Lemma shows that large cliques are not helpful.

\begin{lemma}\label{lemma:cliques} 
Let $p\in\{3,4,5,6,7,8,9\}$ and let $m$ be given. Assume $A\in\Av(m,F(0,p,1,0))$. 
Then  all components of  $G(A)$ induced by the undirected edges are cliques or there is a deletion of cost at least 1. Also $\ncols{A}\le c_pm+1$.
\end{lemma}

\begin{proof} The proof uses induction on $m$. Assume for $m'<m$ and 
$A'\in\Av(m',F(0,p,1,0))$ that the components of $G(A')$ are all cliques or there is a deletion of cost at least 1. Also assume $\ncols{A'}\le c_pm'+1$.

  
  Let $C_1,C_2,\ldots$ be the components of $G(A)$ induced by the undirected edges.  When  each $C_i$ is shrunk to a single vertex then the components $C_i$ can be ordered by the transitive ordering of the resulting directed graph.  Reorder the rows of $A$ to
respect this order. Given two rows $x,y$ with $x\in C_i$
and $y\in C_j$ with $i<j$ in the ordering, then $x\rightarrow y$ hence there is
no submatrix $\linelessfrac{x}{y} \bigl[\linelessfrac{0}{1}\bigr]$.

 For a component $C_i$ on $k$ vertices,  let $B_i$  denote  the $k$-rowed submatrix of
$A_{C_i}$  where $B_i$ is all non-constant columns of 
$A_{C_i}$ plus a column of 0's if present in $A_{C_i}$.
We note that if there is a column non-constant in $B_i$,
then such a column is forced have 1's on rows in components $C_{\ell}$ with $\ell<i$ and 0's on rows in components $C_{\ell}$ with $\ell>i$.  This shows that a column cannot be non-constant on two components. Thus $A$ has the structure of  \rf{construction}. Also it forces $B_i$ to be simple.

 We first show that small components do not occur in $G(A)$ else there is a deletion of cost 1. Since $B_i$ is simple and $k$-rowed, the number of  non constant columns in $B_i$ is at most $2^k-2$.
A component of size $k$ has at most $2^k-2$ non-constant columns.  Thus by
 \rf{costcomponent} a component of size $k$ costs at least $c_p k-(2^k-1)$.   The costs are at least 1 for $p=3,k\le 2$, $p=4,5,k\le 3$ and $p=6,7,8,9$ and $ k\le 4$.

The Upper Bound Lemma shows that there are unique solutions  $K_k\backslash\1_k$ of cost 0 for $p=3$, $k=3$, and $p=5$, $k=4$ and $p=9$, $k=5$ and they  are cliques.  Thus any component of cost 0 is a clique in the cases $p=3$, $k=3$, and $p=5$, $k=4$ and $p=9$, $k=5$.  

For $p=4$, $k=4$, there are exactly  10 columns of sum 1,2,3 with no
$\linelessfrac{1}{2}\left[\linelessfrac{0}{1}\right]$ and the resulting matrix has $F(0,4,1,0)$. Thus for $p=4$,  a component of cost 0 must be a clique.    

For $p=6,7,8$ and $k=5$, a column of sum 2 or 3 on 5 rows produces $6$ configurations $\left[\linelessfrac{0}{1}\right]$  and a column of sum 1 or 4 
produces $4$ configurations $\left[\linelessfrac{0}{1}\right]$.  Consider $p=6$ and $A\in\Av(5,F(0,6,1,0))$. If either $I_5$ or $I_5^c$ in $A$, this yields a clique so assume there are at most 8 columns of sum 1 or 4. There are at most
100 configurations $\left[\linelessfrac{0}{1}\right]$. Let $a$ be the number of columns of sum 1 or 4 and let $b$ be the number of columns of sum 2 or 3.  Then $4a+6b\le 100$. 
Given $a\le 8$ yields $a+b\le 19 $.    This yields  $\ncols{A}\le 21$ and so  yields a component has cost at least 1.    

Now for $p=7,8$, $k=5$,  we can argue as follows.  Assume  \hfil\break$A\in\Av(5,F(0,p,1,0))$ and yet $G(A)$ is not a clique. Assume the directed edge $1\rightarrow 2$ is in $G(A)$. Then the number of columns on 5 rows with no submatrix $\linelessfrac{1}{2}\left[\linelessfrac{0}{1}\right]$ is 24. Given our constructions 
\rf{F7construction} and \rf{F8construction}, we deduce that such a case yields a component of cost at least 1. 

Thus we may assume non-clique components of these small sizes do not occur verifying the Lemma in these cases. Moreover if the cost of such a component is 0 then it is a clique.   

If we have a component
$C_i$ on $k$ vertices with $h$ columns non-constant on $B_i$, then we can
delete the $k$ rows and $\leq h+1$ columns (the $h$ columns non-constant on $B_i$ and possibly 
column of 0's), that are all 1's on components above $C_i$ and all 0's on the 
components below $C_i$,  to obtain a simple $(m-k)$-rowed
matrix.  As noted, in \rf{construction}, $B_i$ consists of non-constant columns plus (typically) the column of 0's and so $h+1$ is (a bound on) the number of columns in $B_i$.  If there is a single component,  then we have the non-constant columns plus  the column of 0's plus the column of 1's and so $h+2$ columns.  Note that $h,h+1,h+2$ are all relevant and this can confuse the reader.

By the Deletion Lemma we are done if
\begin{equation}
    h+2\le  c_p\cdot k,
\label{basicbound}\end{equation}
where the cost of the deletion would be $c_p\cdot k-(h+1)$ and hence at least 1 or, if only one component, we have the correct bound.  
If a component $C$ is a clique, then we can use the Upper Bound Lemma and verify that the cost of the component  is at least 1 except in the small cases noted above ($p=3$, $k=3,4,5$) ($p=4,5$,$k=4$)  and ($p=6,7,8,9$, $k=5$) all of which are cliques. We can construct a table \rf{p=6upperboundlemma} using the Upper Bound Lemma. This is also done in Case 1 below where the Upper Bound Lemma still holds even though $C$ is not a clique.

 The remainder of the proof considers    a component  $C$ (induced by the undirected edges)  that is not a clique. 
  Of course $F(0,p,1,0)\not\prec A_C$ and \rf{wim} applies.
  Recall that $A_C$ does not have a column of 1's but can have the column of 0's.  The argument splits into two cases.

\vskip 10pt
\noindent\textbf{Case 1.} 
$F(0,2p-1,0,0)\not\prec A_C$.

Assume the component has $k$ vertices and consider the $k$-rowed matrix formed from the possible non-constant columns on these $k$ rows. Then  for any pair
of rows $i,j$, to avoid $F(0,p,1,0)$, is either in first case of \rf{wim} and has at most $p-1$ submatrices $\left[\linelessfrac{0}{1}\right]$ and $p-1$ submatrices $\left[\linelessfrac{1}{0}\right]$ or we are in the second or third cases of \rf{wim} with at most $2p-2$ submatrices of one of the two types. Thus rows $i,j$ have  at most $2p-2$ configurations $\left[\linelessfrac{0}{1}\right]$.  The 
Upper Bound Lemma now applies, even though the component need  not be a clique,  with $t=2p-2$.  Thus the maximum number of columns  non-constant on $C$ is
\begin{equation}
    h\le\lfloor 2k+\frac{(2p-6)k(k-1)}{4(k-2)}\rfloor.
\label{basicboundp}\end{equation}
Note that from \rf{basicbound}, $h+1$, $h+2$ are important.   In \rf{construction}, our components $B_i$ are of size $h+1$ while, if there is only one component, the construction yields $h+2$ columns.

We now consider the various values of $p$. Using \rf{basicboundp}  for $p=3$, $k=4$, the component $C$ has $h\le \lfloor 2\cdot 4\rfloor$ and $h+1\le 9$ so the cost is at least $\frac{7}{3}\cdot 4-9=\frac{1}{3}$ which can be achieved ($[K_4^3K_4^1K_4^0]$). Similarly for $p=3$, $k=5$ and the Upper Bound Lemma, component $C$ has $h\le \lfloor 2\cdot 5\rfloor$ and $h+1\le 11$ so the cost is at least $\frac{7}{3}\cdot 5-11=\frac{2}{3}$ which can be achieved ($[K_5^3K_5^1K_5^0]$).  Note that for $p=3$, $k\ge 6$,  the Upper Bound Lemma yields that component $C$ has $h\le \lfloor 2\cdot k\rfloor$ and $h+1\le 2k+1$ so the cost is at least $\frac{7}{3}\cdot k-(2k+1)\ge 1$.
So Case 1 could occur for $k=4,5$ but only with the given  costs if the component is a clique.  One can verify  cost is bigger than one if the component is not a clique.  For $k\ge 6$, the cost is greater than one. The extremal matrices $\ext(m,F(0,3,1,0))$ are described in
\trf{fullF(0,3,1,0)}.

For $p=4$ we can't achieve the bounds of the Upper Bound Lemma  for $k=4,5$.  A computer search was done by a program written by Miguel Raggi \cite{Raggi} written as part of his PhD at UBC in 2011. It has been tested a number of times.  Tackling cases with $m=6$ has not worked for these problems.  We obtain
$\forb(4,F(0,4,1,0))=12$ (cost=0 using $c_4=\frac{11}{4}$),
$\forb(5,F(0,4,1,0))=14$ (cost=$\frac{3}{4}$) and both are only achieved by cliques where other components have cost at least 1.
Now using $p=4$ in \rf{basicboundp}, for $k\ge 6$, the costs for the components are at least 1. So Case 1 does not occur for $p=4$.

 For $p=5$, computer search \cite{Raggi} yields
$\forb(4,F(0,5,1,0))=16$ ($h=14$, cost=0),  
$\forb(5,F(0,5,1,0))=18$ ($h=16$, cost=$1\frac{3}{4}$).  Using \rf{basicboundp} yields the same bounds  for $k=4,5$. 
The cost of  components for $k\ge 6$ is at least 1 by 
\rf{basicboundp}.   Thus Case 1 does not occur for $p=5$.

For $p\in\{6,7,8,9\}$, the cost of a component  using \rf{basicboundp}  is at least
 $$c_p\cdot k-(h+1)= c_p\cdot k- \left(\left\lfloor 2k+\frac{(2p-6)}{4}\frac{k(k-1)}{(k-2)}\right\rfloor+1\right),$$
 which is increasing in $k$ and is minimized at $k=5$.  We have already handled components of size $k=5$ above. 
  The following table indicates that we have cost at least 1 for components of size $k\ge 6$ with the exception of $p=7$ and $k=6$ which requires a more careful analysis.  Recall the cost is $c_p\cdot k-(h+1)$.

\begin{E}
\hskip .7in
\hbox{\hfil
\vbox{\offinterlineskip
\halign{&\vrule#&
\strut\,\,\hfil#\hfil\,\,\cr
\noalign{\hrule}
height2pt&\omit&&\omit&&\omit&&\omit&&\omit&&\omit&\cr
&$p$&&$\,\,c_p\,\,$&&$h$ ($k=6$)&&cost ($k=6$)&&$h$ ($k=7$)&&cost ($k=7$)&\cr
height2pt&\omit&&\omit&&\omit&&\omit&&\omit&&\omit&\cr
\noalign{\hrule}
height2pt&\omit&&\omit&&\omit&&\omit&&\omit&&\omit&\cr
&6&&$\frac{21}{5}$&&23&&$1\frac{1}{5}$&&26&&$2\frac{1}{5}$&\cr
height2pt&\omit&&\omit&&\omit&&\omit&&\omit&&\omit&\cr
\noalign{\hrule}
height2pt&\omit&&\omit&&\omit&&\omit&&\omit&&\omit&\cr
&7&&$\frac{24}{5}$&&27&&$\frac{4}{5}$&&30&&$2\frac{3}{5}$&\cr
height2pt&\omit&&\omit&&\omit&&\omit&&\omit&&\omit&\cr
\noalign{\hrule}
height2pt&\omit&&\omit&&\omit&&\omit&&\omit&&\omit&\cr
&8&&$\frac{27}{5}$&&30&&$1\frac{2}{5}$&&35&&$2\frac{2}{5}$&\cr
height2pt&\omit&&\omit&&\omit&&\omit&&\omit&&\omit&\cr
\noalign{\hrule}
height2pt&\omit&&\omit&&\omit&&\omit&&\omit&&\omit&\cr
&9&&$\frac{31}{5}$&&34&&$2\frac{1}{5}$&&39&&$3\frac{2}{5}$&\cr
height2pt&\omit&&\omit&&\omit&&\omit&&\omit&&\omit&\cr
\noalign{\hrule}
}
\hfil
}  }
\label{p=6upperboundlemma}\end{E}

 For $p=7$ cost of component of size $k$ is
 $\frac{24}{5}\cdot k-\lfloor 2k+\frac{8}{4}\frac{k(k-1)}{k-2}\rfloor-1$
 which is only $\frac{4}{5}$ for $k=6$ but already $2\frac{3}{5}$ for $k=7$.  
We need some special argument for $k=6$. 
 We can have $h=27$ by taking $\left[K_6^5K_6^5K_6^1\right]$. However this is a clique.  On the other hand, if on 6 rows we have a column with 3 1's and 3 0's, then the bound in the Upper Bound Lemma becomes $2\cdot 6 \cdot 5+(n-13)\cdot 8+9\le 6\cdot 6 \cdot 5$, which implies $n\le 26$ and cost $>1$. Similarly, if we only have columns from $\left[K_6^5K_6^4K_6^2K_6^5\right]$, but at least one pair $i,j$  with at most 11 configurations of $\linelessfrac{i}{j}\left[\linelessfrac{0}{1}\right]$, then the bound in the Upper Bound Lemma becomes $2\cdot 6 \cdot 5+(n-12)\cdot 8\le 6\cdot 6 \cdot 5-1$, which implies $n\le 26$ and cost$>1$. Thus, to have 27 non-constant columns on 6 rows that is not a clique, we must have only columns from $\left[K_6^5K_6^4K_6^2K_6^1\right]$ and at least one pair $i,j$  with 12 columns of 
 $\linelessfrac{i}{j}\left[\linelessfrac{1}{0}\right]$. However, for a given pair we only have 10 columns in $\left[K_6^5K_6^4K_6^2K_6^1\right]$ with 
$\linelessfrac{i}{j}[\linelessfrac{1}{0}]$. 
Thus the cost of a non-clique component is at least 1 for $k=6$ and as noted for all $k\ge 7$.    So we can assume Case 1 does not occur for $p=7$. 

This completes Case 1 for $3\le p\le 9$.

\vskip 10pt
\noindent\textbf{Case 2.} $F(0,2p-1,0,0)\prec A_C$.

 Consider
 two rows $i,j\in C$ with $F(0,2p-1,0,0)\prec A_{\{i,j\}}$.   Then we do not have the undirected edge $i-j$.    We may assume we have the submatrix $(2p-1)\cdot\linelessfrac{i}{j}[\linelessfrac{1}{0}]$ and $i\rightarrow j$. 
 For every other row $s$ of $C$ we have either $i\rightarrow s$ or $s\rightarrow j$
since there will either be $p$ 0's in row $s$ in same columns with  $1's$ or yielding
$i\rightarrow s$ (using \rf{wim}) or $p$ 1's in row $s$ in same columns with the  0's yielding
$s\rightarrow j$ (by \rf{wim}). 

 Let the shortest path of undirected edges in $C$ joining $i,j$, denoted  $i=v_1v_2\ldots v_r=j$. Since for all $2\leq a\leq r-1$ either $v_1\rightarrow
v_a$ or $v_a\rightarrow v_r$, then $i-s$, $s-j$ is impossible and hence $r\geq 4$.  
By the minimality of the path $v_1$ and $v_{r-1}$ can't be joined by an undirected edge so either $v_1\rightarrow v_{r-1}$ or
$v_{r-1}\rightarrow v_1$. The latter is impossible by transitivity. Continuing
by induction shows that $v_s\rightarrow v_t$ if $s+1<t$.

Let $R=\{v_1,v_2,\ldots ,v_r\}$ and let $A_R$ denote the submatrix of $A_C$ on rows $R$. 
There are $2r-2$ possible types of non-constant columns and 2 constant columns  on the $r$ rows  
$v_1,\ldots,v_r$ given in \rf{eq:AR-columns}. Remember that $A_C$ may have more than just the $r$ rows and so the columns in $A_R$ have multiplicities in $A_R$ as labelled.

\begin{E}
A_R\quad \begin{array}{c@{}}
\\ v_1\\ v_2 \\v_3\\ v_4\\ \vdots\\v_{r-1}\\ v_r\\ \end{array}
\begin{array}{@{}c}\begin{array}{ccccccccccc}
a_1^1&a_2^0&a_2^1&a_3^0&a_3^1&\cdots&a_{r-1}^0&a_{r-1}^1&a_r^0&b_0&b_1\,\,\\
\end{array}\\
\left[\begin{array}{@{}ccccccccccc}
\,0\,\,&\,1\,\,&\,1\,\,&\,1\,\,&\,1\,\,&\cdots&\,\,1\,\,\,&\,\,1\,\,\,&\,\,1\,\,&\,\,0\,&\,\,1\,\,\\
1&0&0&1&1&\cdots&1&1&1&0&1\\
0&0&1&0&0&\cdots&1&1&1&0&1\\
0&0&0&0&1&\cdots&1&1&1&0&1\\
\vdots&\vdots&\vdots&\vdots&\vdots&&\vdots&\vdots&\vdots&\vdots&\vdots\\
0&0&0&0&0&\cdots&0&0&1&0&1\\
0&0&0&0&0&\cdots&0&1&0&0&1\\
\end{array}\right]\\
\end{array}
\label{eq:AR-columns}\end{E}
Let $a_k^d$ be the number of columns in $A$ with a 0 in row $v_k$ and $d$ in
row $v_{k+1}$ with 1's above and 0's below as in \rf{eq:AR-columns}. We say `columns $a_s^d$\,\,' to refer to these columns.  Using the edges $v_s - v_{s+1}$ we
get that
\begin{align*}
    1\leq a_{s-1}^1+a_{s+1}^0+a_{s+1}^1\leq p-1\;\textrm{for}\;2\leq s\leq r-2\\
    1\leq a_2^0+a_2^1\leq p-1\;\textrm{and}\; 1\leq a_{r-2}^1+a_{r}^0\leq p-1
\end{align*}
Adding these inequalities together gives
\begin{equation*}
    a_1^1+2a_2^1+\ldots+2a_{r-2}^1+a_{r-1}^1+a_2^0+\ldots+a_r^0\leq (p-1)(r-1). 
\label{boundextended}\end{equation*}
Using that $a_s^1\geq 1$ for all $s=1,2,\ldots ,r-1$ (using $v_s - v_{s+1}$) we conclude that
\begin{equation}\label{eq:sumAsdwithp}
    \sum_{s,d} a_s^d=a_1^1+a_2^1+\ldots+a_{r-2}^1+a_{r-1}^1+a_2^0+\ldots+a_r^0\leq (p-2)r-p+4. 
\end{equation}

Among pairs of rows $u,v$ with  $(2p-1)\cdot\linelessfrac{u}{v}\left[\linelessfrac{1}{0}\right]$, the pair $i,j$ is chosen as the pair with the shortest path of undirected edges joining them. Thus  we can't have
$(2p-1)\cdot\linelessfrac{v_1}{v_{r-1}}\left[\linelessfrac{1}{0}\right]$ with $v_1,v_{r-1}$ joined by a path of undirected edges or 
$(2p-1)\cdot\linelessfrac{v_2}{v_{r}}\left[\linelessfrac{1}{0}\right]$ with $v_2,v_{r}$ joined by a path of undirected edges. Now 
$A$ has $t\cdot\linelessfrac{v_1}{v_{r-1}}\left[\linelessfrac{1}{0}\right]$ where
$t= \sum a_s^d-(a_1^1+a_{r}^0)$ and has $t\cdot\linelessfrac{v_2}{v_r}\left[\linelessfrac{1}{0}\right]$ where
$t= \sum a_s^d-(a_2^0+a_2^1+a_{r-1}^1)$. The former yields
$\sum a_s^d-(a_1^1+a_{r}^0)\le 2p-2$. Given $a_1^1\le p-1$ and  $a_{r}^0\le p-1$, we obtain 
\begin{E}\sum a_s^d\le 4p-4.\label{newbound}\end{E}
The latter yields 
$\sum a_s^d-(a_2^0+a_2^1+a_{r-1}^1)\le 2p-2$. Given $a_2^0+a_2^1\le p-1 $ and  $a_{r-1}^1\le p-1$, we again  obtain
\rf{newbound}.  Sharper estimates for $a_1^1+a_r^0$, $a_{r}^0$, $a_{r-1}^1$ and $a_2^0+a_2^1$ can improve the bound \rf{newbound}.  This bound does not depend on $r$ and so is helpful for larger $r$.

There is a decomposition that can assist.
Let $D$ denote the $r\times  \sum_{s,d} a_s^d$ submatrix of the non-constant columns $A_R$ (with multiplicities). From our observation at the beginning of Case 2, 
every  row $s$ of $C\backslash R$ has either $i\rightarrow s$ or $s\rightarrow j$.  Define $V$ as those $s\in C\backslash R$ with  $i\rightarrow s$ and define $U$ as those 
$s\in C\backslash R$ with  $s\rightarrow j$ but with $s\notin V$. Hence $R\cup U\cup V=C$.  Combine that with the $b_0$,$b_1$ columns from \rf{eq:AR-columns} and decompose $A_C$ using rows $R,U,V$ and columns non constant on $R$, columns 0 on rows $R$ and columns 1 on rows $R$. as follows. 
 Row $s\in V$ of columns $b_0$ is all 0's using $s\rightarrow v_r$ and  row $s\in U$ of columns $b_1$ is all 1's using 
 $v_1\rightarrow s$.  This yields the following decomposition of $A_C$ where $\0$ and $\1$ denote blocks of 0's and 1's.

\begin{E}A_C=\begin{array}{c@{}}R\\ U\\ V\\ \end{array}\left[\begin{array}{c|c|c}
 D&\0&\1\\
 \hline
 E&B_0&\1\\
 \hline
 F&\0&B_1\\ 
\end{array}\right],\label{smalldecomp}\end{E}
 We assume $B_1$ does not have the column of 1's so that $A_C$ does not have the column of 1's. In this decomposition let $u$ be the number of rows in $B_0$ and $v$ be the number of rows in $B_1$. Then
$r+u+v= k$.

From \rf {smalldecomp}, $\ncols{A_C}=\ncols{D}+\ncols{B_0}+\ncols{B_1}$. 
From \rf{eq:sumAsdwithp}, 
\begin{E}\ncols{D}=\sum_{s,d} a_s^d\le (p-2)r-p+4.\label{bound}\end{E} 
Since $B_0\in\Av(u,F(0,p,1,0))$, $B_1\in\Av(v,F(0,p,1,0))$ are simple then by induction
$\ncols{B_0}\le\forb(u,F(0,p,1,0))=\lfloor c_pu\rfloor+1$ and $\ncols{B_1}\le\forb(v,F(0,p,1,0))-1=\lfloor c_pv\rfloor$ ($B_1$ doesn't have column of 1's).  

We can delete $|C|$ rows and $\ncols{D}+\ncols{B_0}+\ncols{B_1}$ columns. This deletion has cost at least: 
\begin{E}(r+u+v)c_p-((p-2)r-p+4)- \lfloor c_pu+1\rfloor-\lfloor c_pv\rfloor.\label{newinduct}\end{E}
This is at least 1 for $p\in\{3,4,5,6\}$ and so Case 2 does not occur for $p=3,4,5,6$. 


Deleting $r$ rows $R$ and the appropriate columns to maintain simplicity should have cost at least 2 to yield cost at least 1 for component $C$ using $\rf{newinduct}$. Note the column that is 0's on $U$ and 1's on $V$ could appear twice in $A_{U\cup V}$.
We use both \rf{bound} and  \rf{newbound} in our arguments.  

For $p=7$ with $c_7=\frac{24}{5}$   
then  $\sum a_s^d\le (p-2)r-p+4=5r-3$ and 
$c_pr-((p-2)r-p+4)=\frac{24}{5}r-5r+3=-\frac{1}{5}r+3\ge 2$ for $r\le 5$.  
As noted above,  our argument requires cost at least 2 for $R$.
Using \rf{newbound}, we have $c_pr-(4p-4)=\frac{24}{5}r-24 \ge 2$ for $r\ge 6$.
This concludes the argument for $p=7$.

For $p=8$ with $c_8=\frac{27}{5}$,  
then  $\sum a_s^d\le (p-2)r-p+4=6r-4$ and 
$c_pr-((p-2)r-p+4)=\frac{27}{5}r-6r+4=-\frac{3}{5}r+4\ge 2$  for $r= 4$.   
Using \rf{newbound}, we have $c_pr-(4p-4)=\frac{27}{5}r-28\ge 2$ for $r\ge 6$.  
More detail is required for $r=5$.

$$\begin{array}{cc}
&\begin{array}{@{}cccccccc}
\,\,\,a_1^1&a_2^0&a_2^1&a_3^0&a_3^1&a_4^0&a_4^1&a_5^0\\
\end{array}\\
\begin{array}{c}v_1\\v_2\\v_3\\ v_4\\v_5\\ \end{array}
&\left[\begin{array}{cccccccc}
\,\,0\,&\,1\,&\,1\,&\,1\,&\,\,1\,&\,\,1\,&\,\,1\,&\,\,1\,\\
1&0&0&1&1&1&1&1\\
0&0&1&0&0&1&1&1\\
0&0&0&0&1&0&0&1\\
0&0&0&0&0&0&1&0\\
\end{array}\right]
\end{array}$$

Our choice of $i,j$ for Case 2 yields no $F(0,2p-1,0,0)$ on pairs of rows $v_2,v_5$ and $v_1,v_4$ yielding $a_1^1+a_3^0+a_3^1+a_4^0+a_5^0\le 2p-2$ and
$a_2^0+a_2^1+a_3^0+a_4^0+a_4^1\le 2p-2$. Adding  yields:

\begin{E}\sum a_s^d + (a_3^0+a_4^0)\le 4p-4.\label{newboundk=5}\end{E}

If $a_3^0+a_4^0\ge 3$ then  \rf{newbound}, using \rf{newboundk=5}, reduces to $4p-7=25$ which  handles case $r=5$, $p=8$.   
Now  $a_2^1=a_3^1=1$ is forced else we can reduce \rf{bound} from $6r-4$ to $6r-5$ which  handles case $r=5$ for $p=8$. But if $a_3^0+a_4^0\le 2$ and $a_2^1=a_3^1=1$, we can use an  alternative to our  bounds which is to use the Deletion Lemma on  \rf{eq:AR-columns}. We can delete  row $v_3$ in $A_C$ and at most 4 columns $a_2^1$, $a_3^1$ and $a_3^0+a_4^0$ yielding a deletion of cost at least 1 since  $c_8=\frac{27}{5}\ge 5$.  This concludes the argument for $p=8$.

For $p=9$ with $c_9=\frac{31}{5}$ then by \rf{newbound}, $c_9r-(4p-4)=\frac{31}{5}r-32\ge 2$ for $r\ge 6$. Note the bound is only 1 away from eliminating  $p=9$, $r=5$.   
Now \rf{bound} is $\sum a_s^d \le 7r-5$ and unfortunately
$\frac{31}{5}r-(7r-5)=-\frac{4}{5}r +5< 2$ for $r\ge 4$ but again the bound is only 1 away from eliminating  $p=9$, $r=4,5$.   For $r=5$, we use the above analysis noting that either $a_2^1\ge 2$ or $a_3^1\ge 2$ results in a drop  in the bound to $7r-6$ which would have 
$\frac{31}{5}r-(7r-6)=-\frac{4}{5}r +6\ge 2$ for $r =4,5$. Thus we may assume $a_2^1=a_3^1=1$. If $a_3^0+a_4^0\ge 3$ then  \rf{newbound},  using \rf{newboundk=5},  reduces to $4p-7=29$ which  now handles case $r=5$, $p=9$.  Thus assume  $a_3^0+a_4^0\le 2$, and as above  the deletion of   row $v_3$ in $A_R$ and at most 4 columns $a_2^1$, $a_3^1$ and $a_3^0+a_4^0$ has  cost at least 1 since  $c_9=\frac{31}{5}\ge 5$. 

The case $p=9$, $r=4$ has more work. We would need $\sum a_s^d \le 7r-6=22$.  Because the bound is so close we can assume 
$\ncols{B_0}=\lfloor c_pu+1\rfloor$ and $\ncols{B_1}=\lfloor c_pv\rfloor$ where $B_1$ is missing the column of 1's. Both $B_0$ and $B_1$ are `extremal' and we can use induction (since $u,v<m$) to assert that  the graphs
$G(B_0)$ and $G(B_1)$ have the components induced by undirected edges to be cliques since a
 deletion of cost at least 1 contradicts $\ncols{B_0}=\lfloor c_pu+1\rfloor$ and $\ncols{B_1}=\lfloor c_pv\rfloor$.  Now $B_0$ having all clique components with $\ncols{B_0}=\lfloor c_pu+1\rfloor$ implies, by \rf{extremalclique}, that  $B_0$ must have all cliques of size 5 and moreover have $K_5\backslash \1_5$ on the rows of each component as follows:

$$\begin{array}{c}R\\ \\ U\\  \\  V\\ \end{array}
\left[\begin{array}{c|c|c|c|c}
D&\0&\cdots &\0&\1\\
\hline
*&K_5\backslash \1_5&&&\1\\
\hline
*&&K_5\backslash \1_5&&\1\\
\hline
*&&&K_5\backslash \1_5&\1\\
\hline
*&\0&\cdots&\0&B_1\\ 
\end{array}\right]$$
We note that any row of $K_5\backslash \1_5$ has 15 1's (and 16 0's). 
Thus all entries in rows $U$ in initial columns under $D$ must be 1: If there is a 0 in row $s\in U$ then there is a 1 in that column in $D$ in some row  $t\in R$ (e.g. in rows $v_1$ or $v_2$)  and this yields $F(0,9,1,0)$ in rows $s,t$.  The same argument would apply to $B_1$ and rows of $V$, interchanging roles of 0's and 1's.   Thus all entries in rows $V$ in initial columns under $D$ must be  all 0's. But then all column multiplicities in $A_R$ for the non constant columns would be at most 1.  This yields $\sum a_s^d=2r-2=8$ and so the cost deleting rows $R$ is at least 2.

Thus Case 1 and Case 2 do not occur and hence components either are cliques or there is a deletion of  cost at least 1.

We now verify  $\ncols{A}\le c_pm+1$.  If all components are cliques then \lrf{cliquesgivebound} gives the bound.  If we have a deletion of cost at least 1 then we can delete $k_1$ rows and $t_1$ columns  with $c_pk_1-t_1\ge 1$.  After the deletion  we  obtain an $A^{(1)}\in\Av(m-k_1,F(0,p,1,0))$ with $\ncols{A^{(1)}}=\ncols{A}-t_1$.   We apply induction with $m-k_1<m$ to deduce that  $\ncols{A^{(1)}}\le c_p(m-k_1)+1$ and this yields $\ncols{A}\le c_pm<c_pm+1$.
 \qed  \end{proof}

\vskip 20pt

\section{ Extremal Matrices for $F(0,3,1,0)$ }\label{F2(0,3,1,0)}

We give  a detailed description for $\ext(m,F(0,3,1,0))$.  
\begin{thm}Assume $m\ge 3$. Then $\forb(m,F(0,3,1,0))=\lfloor\frac{7}{3}m\rfloor+1$.
\label{p=3}\label{fullF(0,3,1,0)}\end{thm}

\proof The upper bound is established using \lrf{lemma:cliques} in combination with \lrf{risefall}.
We have the base cases $\forb(3,F(0,3,1,0))=8$ with construction $K_3$, $\forb(4,F(0,3,1,0))=10$ with construction $[K_4^0K_4^1K_4^3K_4^4]$ and $\forb(5,F(0,3,1,0))=12$ with construction $[K_5^0K_5^1K_5^4K_5^5]$.  The constructions are forced to be unique by Upper Bound Lemma for cliques.   For $m\equiv 0(\hbox{mod }3)$,  take all components to be cliques of size 3 with $B_i=[K_3^0K_3^1K_3^2]$ in \rf{construction}. For  $m\equiv 1(\hbox{mod }3)$ use the construction $[K_4^0K_4^1K_4^3]$  on $4$ rows with all remaining components to be cliques of size 3.  For  $m\equiv 2(\hbox{mod }3)$ use the construction $[K_5^0K_5^1K_5^4]$  on $5$ rows with all remaining components to be cliques of size 3.  Thus 
$\forb(m,F(0,3,1,0))\ge\lfloor\frac{7}{3}m\rfloor+1$ for $m\ge 3$.
\qed\vskip 10pt

\begin{thm} Let $A\in\ext(m,F(0,3,1,0))$ using \rf{construction} for notation. For $m\equiv 0(\hbox{mod }3)$, the components of $G(A)$ are cliques of size 3 . For $m\equiv 1(\hbox{mod }3)$, the components of $G(A)$ consist of one clique of size 4 and the rest are cliques of size 3. For $m\equiv 2(\hbox{mod }3)$, the components of $G(A)$ consist of one clique of size 5 and the rest are cliques of size 3 or two cliques of size 4 and the rest are cliques of size 3. 
A clique of size 3 has $B_i=[K_3\backslash \1_3]$, a clique of size 4 has $B_i=[K_4^0K_4^1K_4^3K_4^4]$, and a clique of size 5 has $B_i=[K_5^0K_5^1K_5^4]$. \label{F(0,3,1,0)ext}\end{thm}

\begin{proof} By \trf{fullF(0,3,1,0)}, $\forb(m,F(0,3,1,0))=\lfloor\frac{7}{3}m\rfloor+1$. We use the cost idea. The only possible component sizes are $3,4,5$; all other components have cost  at least 1 and hence do not occur.   The cost of a component of size $3$ is 0, the cost of a component of size $4$ is 1/3 and the cost of a component of size 5 is 2/3.

We deduce that for $m\equiv 0(\hbox{mod }3)$, that each component must have cost 0 and hence all components of size 3.   We deduce  that for $m\equiv 1(\hbox{mod }3)$, the sum of the costs of the components is 1/3 and hence there is one component of size 4 and the rest of size 3.    We deduce  that for $m\equiv 2(\hbox{mod }3)$, the sum of the costs of the components is 2/3 and hence there is one component of size 5 and the rest of size 3 or there are two components of size 4 (total cost 2/3) and the rest of size 3. Above are given the unique constructions for each of the component sizes.   
 The different ordering of the components means these examples yield many extremal matrices.  \qed\end{proof}

 \section{ Extremal matrices for $F(0,6,1,0)$ }\label{F2(0,6,1,0)}

 In determining extremal matrices the following is essential. 
\begin{lemma}Let $p\in\{6,7,8,9\}$ and $A\in\ext(m,F(0,p,1,0))$.  Assume all components are cliques. Then either there is a clique whose deletion costs at least 1 or all cliques have size 5. \label{extremalclique}\end{lemma}

\proof Clique components  of size less than 5 or greater than 5  all have cost at least 1 by the Upper Bound Lemma.  
\qed\vskip 10pt

 \vskip 10pt
\noindent{\bf Proof of \trf{F61}}. 
For $p=6$, a clique component of size $k\ne 5$  has cost at least 1 in these cases yielding the bound of \trf{F61}.  
If a component of $k\geq 6$ vertices is a clique then we can apply the Upper Bound Lemma with $t=10$ since there are at most $t=5$
configurations $\linelessfrac{i}{j}\left[\linelessfrac{0}{1}\right]$ in any pair of rows $i,j$ and so at most 10 configurations $\left[\linelessfrac{0}{1}\right]$ in rows $i,j$ ($5$ in each of the two orderings).  
The maximum number of non-constant columns on the clique
is $\lfloor 2k+\frac{6k(k-1)}{4(k-2)}\rfloor$ by \rf{risefallsimple}. Using
\begin{equation*}
    \displaystyle\max_{k:\,k\geq 6}\left\lfloor
2k+\frac{6k(k-1)}{4(k-2)}+1\right\rfloor< \frac{21}{5}\cdot k,
\end{equation*}
combined with \rf{basicbound}, establishes the bound. 
\qed

\vskip 10pt
 All component sizes other than 5 have a cost of at least 1. Note for  $m\equiv 0(\hbox{mod }5)$, a matrix  $A\in\ext(m,F(0,6,1,0))$ must have all components of size 5  and chosen from one of the three matrices in
 $\ext(5,F(0,6,1,0))$ in order for $\ncols{A}=\lfloor \frac{21}{5}m\rfloor+1$ using the construction  \rf{construction} with each $B_i$ being a matrix in $\ext(5,F(0,6,1,0))$ minus the column of 1's.    This also yields that for 
  $m\not\equiv 0(\hbox{mod }5)$ that $\forb(m,F(0,6,1,0))<  
 \lfloor \frac{21}{5}m\rfloor+1$.   
 
 For $k=6$ the cost is $\ge\frac{6}{5}$. 
 A construction for $m=6$ with $\lfloor\frac{21}{5}\cdot 6\rfloor=25$ columns is
 \begin{E}\left[K_6^6K_6^5\,
 \begin{array}{ccccccccccc}
 1&0&1&1&1&0&0&0&0&0&0\\
 1&0&0&0&0&1&1&1&0&0&0\\
 1&0&0&0&0&0&0&0&1&1&1\\
 0&1&1&0&0&1&0&0&1&0&0\\
 0&1&0&1&0&0&1&0&0&1&0\\
 0&1&0&0&1&0&0&1&0&0&1\\
 \end{array}\,K_6^1K_6^0\right].\label{6,F(0,6,1,0)}\end{E}
 \noindent For  $m\equiv 1(\hbox{mod }5)$ and $m> 6$, the construction would have one component of size 6 (remove the column of 1's above)  and the rest of size 5 from
 $\ext(5,F(0,6,1,0))$. \qed
\vskip 10pt
 While determining  $\ext(m,F(0,6,1,0))$ for all $m\ge 7$ might be nice in analogy to \trf{F(0,3,1,0)ext}, we would need $\forb(7,F(0,6,1,0))$, $\forb(8,F(0,6,1,0))$ and $\forb(9,F(0,6,1,0))$  to establish $\forb(m,F(0,6,1,0))$ exactly.  This  computation appears daunting.

 Theorems~\ref{F71}, \ref{F81} and \ref{F91}  for $p=7,8,9$ follow from \lrf{lemma:cliques}.  Note that for $p=9$, we use \lrf{cliquesgivebound} and the only clique component of 5 rows of cost 0 is  $K_5\backslash \1_5$.   Hence $\ext(m,F(0,9,1,0))$ is unique for 
 $m\equiv 0(\hbox{mod }5)$:

$$
\left[\begin{array}{c|c|c|c|c}
K_5\backslash \1_5&\1&\1&\1&\1_5\\
\hline
\0&K_5\backslash \1_5&\1&\1&\1_5\\
\hline
\vdots&&\ddots&&\vdots\\
\hline
\0&\0&\0&K_5\backslash \1_5&\1_5\\
\end{array}
\right]$$

\end{document}